\documentclass[oneside,a4paper,english,links]{amca}
\usepackage{graphicx}
\usepackage{amsmath,amsfonts}
\usepackage[utf8]{inputenc}

\usepackage{epstopdf}

\newtheorem{theorem}{Theorem}[section]
\newtheorem{lemma}[theorem]{Lemma}

\newtheorem{example}[theorem]{Example}
\newcommand{\R}{\mathbb{R}}
\newcommand{\C}{\mathbb{C}}

\newcommand{\disp}{\displaystyle}

\title{A REGULARIZATION OPERATOR FOR THE SOURCE APPROXIMATION OF A TRANSPORT EQUATION}

\author[a,b]{Guillermo F. Umbricht}
\author[b]{Diana Rubio}
\author[a]{Claudio D. El Hasi}
\affil[a]{Instituto de Ciencias, Universidad Nacional de General Sarmiento, Juan María Gutiérrez 1150, 1613 Los Polvorines, Buenos Aires, Argentina, guilleungs@yahoo.com.ar, http://www.ungs.edu.ar}
\affil[b]{Centro de Matemática Aplicada, Escuela de Ciencia y Tegnología. Universidad Nacional de San Martin,25 de mayo y Francia,1650 San Martin, Buenos Aires, Argentina, drubio@unsam.edu.ar, http://www.unsam.edu.ar}

\begin{document}
\vspace{3cm}

\maketitle



\begin{abstract}
Source identification problems have multiple applications in engineering such as the identification of fissures in materials,
 determination of sources in electromagnetic fields or geophysical applications, detection of contaminant sources, among others.
 In this work we are concerned with the determination of a time-dependent source in a transport equation from noisy data measured  at a fixed position. By means of Fourier techniques can be shown that the problem is ill-posed in the sense that the solution exists but it does not vary continuously with the data. A number of different techniques were developed by other authors to approximate the solution.  In this work, we consider a family of parametric regularization operators to deal with the ill-posedness  of the problem. We proposed a manner to select the regularization parameter as a function of noise level in data in order to obtain a regularized solution that approximate the unknown source.  We find a Hölder type bound for the error of the approximated source 
when the unknown function is considered to be bounded in a given norm. Numerical examples illustrate the convergence and stability of the method.
\end{abstract}

\section{INTRODUCTION}

Source identification problems are inverse problems of great interest due to the amount of applications in different disciplines. For instance, in the problems of source identification  we can find applications in heat conduction processes \cite{HOL93}, contaminant detection \cite{Letall2006} and detection of tumor cells \cite{Macleod99}. The determination of the foci of pollution of the layers of groundwater triggers environmental and health problems for the population. The identification of sources of pollution in groundwater can be modeled by a transport equation where its input represents the mean concentration of contaminants per unit average of effective porosity \cite{Letall2006, Sun96}. 
As a first approach, this work focus on the determination of a source for a one-dimensional equation based on noisy measurements taken in a fixed position, in an unbounded domain.
This is an ill-posed problem in the sense of Hadamard \cite{Hadamard23} because its solution does not depend continuously on data. In particular, the high frequency components in arbitrarily small data errors can lead to arbitrarily large errors in the result.
For this type of problems, regularization methods are widely used in the literature in order to obtain a stabilized approximate solution \cite{Fu04,JM08}. 
Here, we designed  a (parametric) family of regularization operators (see \cite{Engel1996}) that compensates the factor that causes the instability of the inverse operator for this problem. This family of operators leads to  well-posed problems that approximates the given ill-posed problem. One must select an {\it a priori} or {\it a-posteriori} rule to choose the regularization parameter, ({\it parameter choice rule}). For more details see \cite{Engel1996}, \cite{Kirsch2011}.
The proposed family of regularization operators provides a framework in the theory of operators for the modified regularization method used, for instance, in  \cite{QianFu06, QianFu07, Xiao11, Zahoelatl14}.\\
Assuming that the source is bounded in a given Hilbert space $H^p$, we propose a parameter choice rule that depends on $p$ and on the data noise level 
.\\
We demonstrate the stability and convergence of the regularization family and obtain a H\"older type bound for the estimation error. 
Numerical examples illustrate its performance including the calculation of the error and the theoretical bound.

\section{The problem of the source determination}

We focus on the  problem of determining the source $f$ for the one-dimensional transport equation from noisy data measurements with conditions. Specifically, we look for the source $f$ that satisfies the system 
\begin{equation}
\label{transpeqn}
\begin{cases} 
u_t(x,t)=\alpha^2 u_{xx}(x,t)-\beta u_x(x,t)-\nu u(x,t)+f(t), \qquad & x \in \R, \quad t>0 \\
u(0,t)=0, \qquad & t>0\\
u(x_0,t)=y(t), \qquad &  t>0, \quad x_0>0 \\ 
\displaystyle \lim_{x\rightarrow\infty} |u(x,t)| < M \qquad & t>0, \mbox{ for some } M\in \R \\
\end{cases}
\end{equation}
where $\alpha^2, \nu >0$, \, $\beta \geq 0$. In addition, we assume that $u(x,\cdot), f(\cdot) \in L^{2}(\R)$ are unknown functions and that $y \in L^2(\R)$ can be measured with certain noise level $\delta$, i.e., the data function $y_{\delta} \in L^{2}(\R)$ satisfies $||y-y_{\delta }||_{L^{2}(\R)}\le \delta.$ 

By means of the Fourier transform, the problem can be written in the frequency space as
\begin{equation}
\label{eqnfreq}
\begin{cases}
-\alpha^2\hat{u}_{xx}(x ,\xi)+\beta\hat{u}_{x}(x ,\xi)+ z(\xi) \, \hat{u}(x,\xi)=\hat{f}(\xi),\qquad  & \xi,x \in \R,  \\
 \hat{u}(0,\xi)=0, & \xi \in \R \\
\hat{u}(x_0,\xi)=\hat{y}(\xi ), & \xi \in \R, \, x_0>0\\
\displaystyle \lim_{x\rightarrow\infty} |\hat{u}(x,\xi)| < M^* \qquad & \xi \in \R, \mbox{ for some } M^*\in \R, \\
\end{cases}
\end{equation}
 where $z(\xi)=\nu + i\xi  \in \C$. 
 
Solving the second order linear ordinary differential equation \eqref{eqnfreq} we obtain
\begin{equation}
\label{illf}
\hat{f}(\xi )=\Lambda (\xi )\hat{y}(\xi ) \qquad  \mbox{ where } \quad \Lambda (\xi )=\frac{z(\xi)}{1-e^{\frac{\beta-\sqrt{\beta^2+4\alpha^2 z(\xi)}}{2\alpha^2 }\, x_0} }.
\end{equation}

We introduce the operator $T$  for $y \in L^2(\R)$,
\begin{equation}
\label{operatorT}
Ty=\frac{1}{\sqrt{2\pi}}\int_{-\infty}^{+\infty}{e^{i\xi t}\Lambda (\xi) \, \hat{y}(\xi)d\xi}. 
\end{equation}
Note that
\begin{equation}
\label{ip}
f(t) =T y (t).
\end{equation}
For simplicity, we drop the variable $t$ in \eqref{ip} and denote $f_\delta = T y_\delta$. Hence
\begin{equation}
  f -f_\delta= T(y -y_\delta )  = \frac{1}{\sqrt{2\pi}}\int_{-\infty}^{+\infty}{e^{i\xi t}\Lambda (\xi) \, (\hat{y}(\xi) - \hat{y}_\delta (\xi)) d\xi}.
\end{equation}


We observe that the factor $\Lambda (\xi)$ increases without bound as $|\xi| \to \infty$ amplifying the high frequency components of the error $ \hat{y}(\xi) - \hat{y}_\delta (\xi)$. This implies that the solution does not depend continuously on the data and the problem is ill-posed in the sense of Hadamard \cite{Hadamard23}. 

Regularization methods are commonly used in dealing with unstable solutions. Regarding the determination of a  source for a parabolic equation, most of the research articles  that use certain regularization techniques are restricted to particular  cases. In \cite{Sivergina2003} the author focus on a convection-diffusion equation, while in  \cite{DF2009} , \cite{DFY2009} , \cite{TLA2005} , \cite{TQA2006} , \cite{YF10} , \cite{Yaetal2010} , \cite{Zahoelatl14} only diffusion is considered. In this work we focus on  a general one dimensional advection-difussion transport equation with a time dependent source.

We define a family of operators and an {\it a-priori}  parameter choice rule \cite{Engel1996} , \cite{Kirsch2011} in order to  regularize the solution to the inverse source problem \eqref{ip}. The stability and convergence of the proposed  regularization family is analyzed and an error bound is obtained based on the data noise level, assuming some smoothness on the source.
The performance of this approach is numerically illustrated. A comparision with the unregularized solution is included.

\section{Inverse problem regularization}

\label{reg}

\smallskip
 
{\bf Definition} \cite{Kirsch2011}
Let $T : Y \longrightarrow X$,  $X$ and $Y$ be Hilbert spaces and $T$ an unbounded operator. A regularization strategy for $T$ is a family of linear and bounded operators 
\begin{equation}
\label{DefinitionR}
R_{\mu} := Y \longrightarrow X, \quad \mu>0, \quad /  \quad  \lim _{\mu \to 0} R_{\mu} y = Ty, \quad \forall y \in Y.
\end{equation}

Let us define the parametric family of integral operators $R_\mu: L^{2}(\R) \to L^{2}(\R)$ for $\mu \in (0,1)$,
\begin{equation}
\label{familyR}
R_{\mu} y := \frac{1}{\sqrt{2\pi}} \int_{-\infty}^{+\infty}{e^{i\xi t} \frac{\Lambda(\xi)}{1+\mu^2\xi^2} \hat{y}(\xi) d\xi}
\end{equation}
where $\Lambda (\xi )$ given in \eqref{illf} and the denominator $1+\mu^2 \xi^2$ is introduced for stabilization purposes. 

\begin{theorem}
Let us consider the problem of identifying $f$ from noisy data $ y_\delta(t)$  measured at a given position $x_0>0$, where  $\delta$ is the data noise level and $u$ and $f$ satisfy
\begin{equation}
\label{illpp}
\begin{cases} 
u_t(x,t)=\alpha^2 u_{xx}(x,t)-\beta u_x(x,t)-\nu u(x,t)+f(t), \qquad & x \in \R, \, t>0 \\
u(0,t)=0, \qquad & x \in \R \\ 
\displaystyle \lim_{x\rightarrow\infty} |u(x,t)| < M \qquad & t>0, \mbox{ for some } M\in \R. \\
\end{cases}
\end{equation}
Let $\{R_{\mu}\}$ be the family of operators defined by \eqref{familyR}.
Then, for every $y(t)=u(x_0,t)$ there exists an a-priori parameter choice rule for $\mu >0$ such that the pair $(R_{\mu}, \mu)$ is a convergent regularization method for solving the identification problem \eqref{ip}.\\
\end{theorem}

{\bf Proof }
The factor $\disp \frac{\Lambda (\xi)}{1+\mu^2 \xi^2}$ is bounded for all $\xi$  since it is continuous for all $\xi \in \R$ and
\[\lim _{\xi \to \pm \infty} \left\vert \frac{\Lambda (\xi)}{1+\mu^2\xi^2} \right\vert 
=\lim _{\xi \to \pm \infty} \left\vert \frac{ z(\xi)}{\left(1-e^{\frac{\beta-\sqrt{\beta^2+4\alpha^2 z(\xi)}}{2\alpha^2 }\, x_0}\right)(1+\mu^2\xi^2)} \right\vert
<\infty. \]
Hence, for all $\mu>0$, $R_{\mu}$ is a continuous operator and $R_{\mu}  \to T$ pointwise on $L^{2}(\R)$ as  $\mu \to 0$  for $T$ given in (\ref{operatorT}). Therefore, by Proposition 3.4 in \cite{Engel1996}, $R_{\mu}$ is a regularization for $T$ and for every $y(t)=u(x_0,t)$ there exists an a-priori parameter choice rule $\mu$ such that $(R_{\mu}, \mu)$ is a convergent regularization strategy for solving \eqref{ip}. The regularized solution is given by
\begin{equation}
\label{ffreg}
f_{\delta ,\mu} = R_{\mu} y_{\delta} =\frac{1}{\sqrt{2\pi}} \int_{-\infty}^{+\infty}{e^{i\xi t} \frac{\Lambda(\xi)}{1+\mu^2\xi^2} \hat{y}_{\delta}(\xi) d\xi}.\\
\end{equation}

\smallskip
 
{\bf Remark  } We observe that the proposed family of regularized operators  \eqref{ffreg} is equivalent to the regularization method proposed in \cite{YF11} when taking $\alpha^2=1$; $\beta =\nu =0$ and $x_0=1$. That is, for  the problem of the source estimation in the heat equation for an isolated bar, the operator  $\Lambda (\xi )=\frac{i\xi}{1-e^{-\sqrt{i\xi}}}$  and the  regularized solution is the same as in \cite{YF11}. \\

\section{Error analysis}
\label{erroranalysis}


In this section we are concerned with the stability and error analysis of the regularization method. 
We assume that the source $f$ is bounded in the Sobolev space $H^{p}(\R), \, p>0$, i.e., for some $C>0$ it holds
\begin{equation}
\label{Boundf}
\|  f\|  _{H^{p}(\R)} := \left(\int_{-\infty}^{\infty }{| \hat{f}|^2 \left(1+\xi^2 \right)^{p}d\xi } \right)^{1/2} 
\le C.
\end{equation}
Optimal order strategies require a parameter choice rule that depends on the  a-priori bound $C$ in \eqref{Boundf}. 
Since in practice this bound is unknown, we look for  a regularization parameter that leads to a convergent strategy and that only depends on the noise level of the data $\delta$.
Some results are now introduced that will be used later to obtain a bound for the regularization error, that is, the error between the source $f(x)$ and its estimate $f_{\delta ,\mu }(x)$. 
\begin{lemma}
\label{lemma1}
For $\omega \in \C$ with $Re(\omega )>0$ holds \\
\hspace*{2cm} $\disp \left\vert \frac{1}{1-e^{-\omega }} \right\vert \le \frac{1}{1-e^{-Re(\omega )}}$ , 
\qquad 
$Re(\sqrt {\omega})=\sqrt {\frac{Re(\omega)+\left\vert\omega    \right\vert}{2}} \geq \sqrt {\ Re {(\omega)}}.$
\end{lemma}
\begin{lemma} 
\label{lemma2}
 If  \quad $0<\mu <1$  then  $\disp \frac{\vert x \vert}{1+x^2\mu^2} \le \frac{1}{2\mu}, \quad  \forall x \in \R$.
 
\end{lemma}
\begin{lemma} 
\label{lemma3}
The function $g:\R_{>0} \to \R$ given by 
$\disp g(x)=
\begin{cases}
  \disp \frac{x}{1-e^{-x}} & 0<x<1\\
  \disp \frac{1}{1-e^{-x}} & 1\leq x
  \end{cases}$
satisfies $g(x) \leq 2 , \,\forall x>0$.
\end{lemma}
\begin{lemma} 
\label{lemma4}
Let $\alpha^2,\nu,x_0 >0$, \, $\beta \geq 0$  and $0<\mu <1$ then, 
\begin{equation*}
 \left\vert  \frac{\Lambda(\xi) }{1+\xi^2\mu^2} \right\vert \le   \frac{1}{\mu^2} \,\left( \frac{2\alpha^2 (2 \nu +1)}{(-\beta+\sqrt{\beta^2+4\alpha^2\nu}) \, x_0}\right).
\end{equation*}
\end{lemma} 
{\bf Proof}
From equation (\ref{illf}), Lemma \ref{lemma1} and the triangular inequality we have   
\begin{equation}
\label{D_ineq}
	\left\vert \frac{\Lambda (\xi)}{1+\xi^2\mu^{2}}\right\vert \le 
\frac{\nu +\vert \xi \vert}{\left(1-e^{-(\frac{-\beta+\sqrt{\beta^2+4\alpha^2\nu}}{2\alpha^2} ) \, x_0}\right)(1+\xi^2\mu^2)}.
\end{equation}
Let us denote $m=m(\alpha,\beta,\nu) = \frac{-\beta+\sqrt{\beta^2+4\alpha^2\nu}}{2\alpha^2}$ and consider two cases: $m \, x_0\ge 1$ and $m\, x_0 \in (0,1)$\\
 
\setlength{\leftskip}{0pt} 
\setlength{\leftskip}{0pt}
 
\emph{Case $m\, x_0 
\ge 1$} : Using Lemmas \ref{lemma2}-\ref{lemma3}, for $0< \mu <1$ results  

\setlength{\leftskip}{2cm}

\begin{align}
\label{part1}
\phantom{space} \frac{\nu +\vert \xi \vert}{\left(1-e^{-m
\, x_0} \right)(1+\xi^2\mu^2)} &\leq 2 \left( \frac{\nu }{1 + \xi^2 \mu^2} + \frac{\vert \xi \vert}{1 + \xi^2 \mu^2} \right) 
\leq 2\nu +\frac{1}{\mu}.
\end{align}

\setlength{\leftskip}{0pt} 
\setlength{\leftskip}{0pt}
 \vspace{.5cm}
 
\emph{Case $m\, x_0 
\in (0,1)$}: \\

\setlength{\leftskip}{1cm}

\noindent Observe that multiplying and diving by $m\, x_0 $
for $0< \mu <1$, Lemmas \ref{lemma2}-\ref{lemma3} imply
\begin{align}
\label{part2}
\frac{\nu+\vert \xi \vert }{\left(1-e^{-m\,  
x_0}\right)(1+\xi^2\mu^2)}
& \le  \frac{2}{m\, x_0} \, \left( \frac{\nu+\vert \xi\vert}{1+\xi^2\mu^2}\right) \, 
\nonumber\\
&  \le   \frac{2}{m\, x_0} \, \left( \frac{\nu}{1+\xi^2\mu^2} + \frac{1}{2\mu}\right) \, 
  \le  \left(2 \nu + \frac{1}{\mu}\right) \,
 \frac{1}{m\, x_0}.
\end{align}
Note that for $0 < \mu <1$,  equations (\ref{part1}) - (\ref{part2}) yield
\begin{equation*}
 \left\vert  \frac{\Lambda(\xi) }{1+\xi^2\mu^2} \right\vert \le   \frac{1}{\mu^2} \,\left( \frac{2\alpha^2 (2 \nu +1)}{(-\beta+\sqrt{\beta^2+4\alpha^2\nu}) \, x_0}\right).
\end{equation*}
and the proof is completed.

\setlength{\leftskip}{0pt}
\begin{theorem}
\label{boundestimate}
Consider the inverse problem of determining the source $f(t)$ in (\ref{transpeqn}). Let $f_{\delta ,\mu }(t)$ be the regularization solution  given in (\ref{ffreg}) and assume that $\|  f \| _{H^p(\R)} $ is bounded in $H^p(\R)$ (\ref{Boundf}).
Then choosing the regularization parameter $\mu^2=\delta^{\frac{2}{p+2}}$,  there exists a constant $K$  independent of $\delta$ such that
\begin{equation}
\label{errorbound}
\|  f -f_{\delta ,\mu } \| _{L^2(\R)} \le K \, \max \left\{ \delta^{\frac{p}{p+2}}, \delta^{\frac{2}{p+2}}\right\}. 
\end{equation}
\end{theorem}
{\bf Proof}
From now on, let us denote $\|   \cdot \|   = \|   \cdot \|  _{L^2(\R)}$. 
Defining $\disp \hat{f}_{\mu}(\xi) :=  \frac{\Lambda(\xi)}{1+\mu^2\xi^2} \hat{y}(\xi) $, one has
\begin{eqnarray}
\label{ineqteo2}
\left| \hat{f}(\xi )-\hat{f}_{\mu}(\xi )   \right|  
&=& \left| \hat{f}(\xi )-\frac{\Lambda (\xi)}{1+\mu^2\xi^2} \hat{y}(\xi )   \right| 
=  \left|    \hat{f}(\xi )\frac{(1+\xi^2)^{\frac{p}{2}}}{(1+\xi^{2})^{\frac{p}{2}}} 
                            \left(1-\frac{1}{1+\xi^2\mu^2}\right) \right| \nonumber\\
& \le & \sup_{\xi \in \R} \left| (1+\xi^{2})^{-\frac{p}{2}}\left(1-\frac{1}{1+\xi^2\mu^2}\right)  
          \right| 		\left|  \hat{f}(\xi) (1+\xi^2)^{\frac{p}{2}} \right|.
\end{eqnarray}
By \eqref{ffreg} and the triangle inequality we have that
\begin{equation}
\label{dnormas}
 \|   \hat{f} - \hat{f}_{\delta ,\mu }  \|   \leq \|   \hat{f} - \hat{f}_{\mu }  \|    +  \|  \hat{f}_{ \mu }   - \hat{f}_{\delta ,\mu }\|.  
\end{equation}
Now, (\ref{ineqteo2})- (\ref{dnormas}), the definition of $H^p(\R)$-norm given in (\ref{Boundf}) and \eqref{ffreg} lead to
\begin{eqnarray*}
\left\|   \hat{f}  -\hat{f}_{\delta ,\mu } \right\| 
 \le  \sup_{\xi \in \R} \left| (1+\xi^{2})^{-\frac{p}{2}}\left(1-\frac{1}{1+\xi^2\mu^2}\right)  
          \right| 	 \|  f \| _{H^p(\R)}
+ \sup_{\xi \in \R} \left|  \frac{\Lambda (\xi)}{1+\xi^2\mu^2} \right|
\left\|    \hat{y}-\hat{y}_{\delta } \right\|. 
\end{eqnarray*}
From \cite{YF10}, it holds 
\begin{equation*} 
\disp \sup_{\xi \varepsilon \R} \left\vert (1+\xi^2)^{-\frac{p}{2}} \left(1-\frac{1}{1+\xi^2 \mu^2} \right)\right\vert \le \max \left\{\mu^p,\mu^2\right\}.
\end{equation*}
Thus, Lemma \ref{lemma4} and  the assumption  $\left\|    \hat{y}-\hat{y}_{\delta }\right\| \leq \delta$ yields to
\begin{equation*}
\|   \hat{f}  -\hat{f}_{\delta ,\mu }  \| 
\le   \max \left\{ \mu ^{p},\mu^2 \right\} 
       \|  f \| _{H^p(\R)} +  \frac{\delta}{\mu^2} \left( \frac{2\alpha^2 (2 \nu +1)}{(-\beta+\sqrt{\beta^2+4\alpha^2\nu}) \, x_0}\right) . 
\end{equation*}
By Parseval's identity, the linearity of the Fourier transform and \eqref{Boundf}, choosing  $\mu^2=\delta^{\frac{2}{p+2}}$  we obtain
\begin{equation}
\label{cota}
 \|  f -R_{\mu } y_{\delta} \| =\left\|   \hat{f}  -\hat{f}_{\delta ,\mu } \right\| \le K  \max \left\{ \delta^{\frac{p}{p+2}}, \delta^{\frac{2}{p+2}}\right\}.
 \end{equation}
 where $K= C +  \left( \frac{2\alpha^2 (2 \nu +1)}{(-\beta+\sqrt{\beta^2+4\alpha^2\nu}) \, x_0}\right)$ and $C$ is the bound for the $H^p(\R)$-norm of $f$, that is, $\|  f\|  _{H^{p}(\R)} \le C$.\\
 
 \smallskip
 
{\bf Remark  }
We notice that although we assumed that $\|  f \| _{H^p(\R)} $ is bounded, in practice we do not know  an {\it a-priori} value for a  bound. For this reason we consider a regularization parameter that does not depends on this bound, it only  depends on the noise data level $\delta$ and $p$.

\section{Numerical examples}
\label{numex}
\begin{figure}[h!]
\includegraphics[scale=0.55]{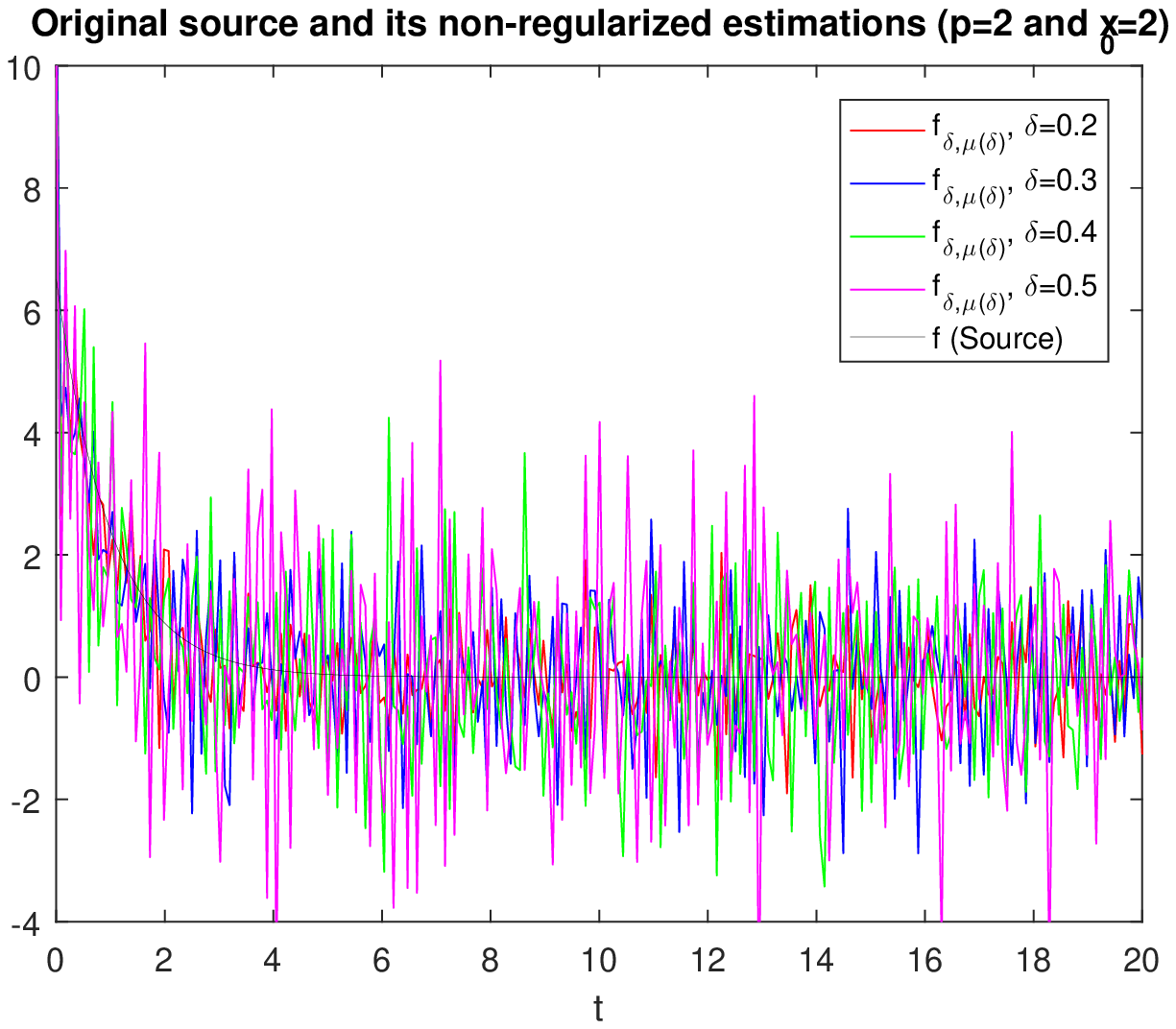}\hspace{-.25cm}
\includegraphics[scale=0.55]{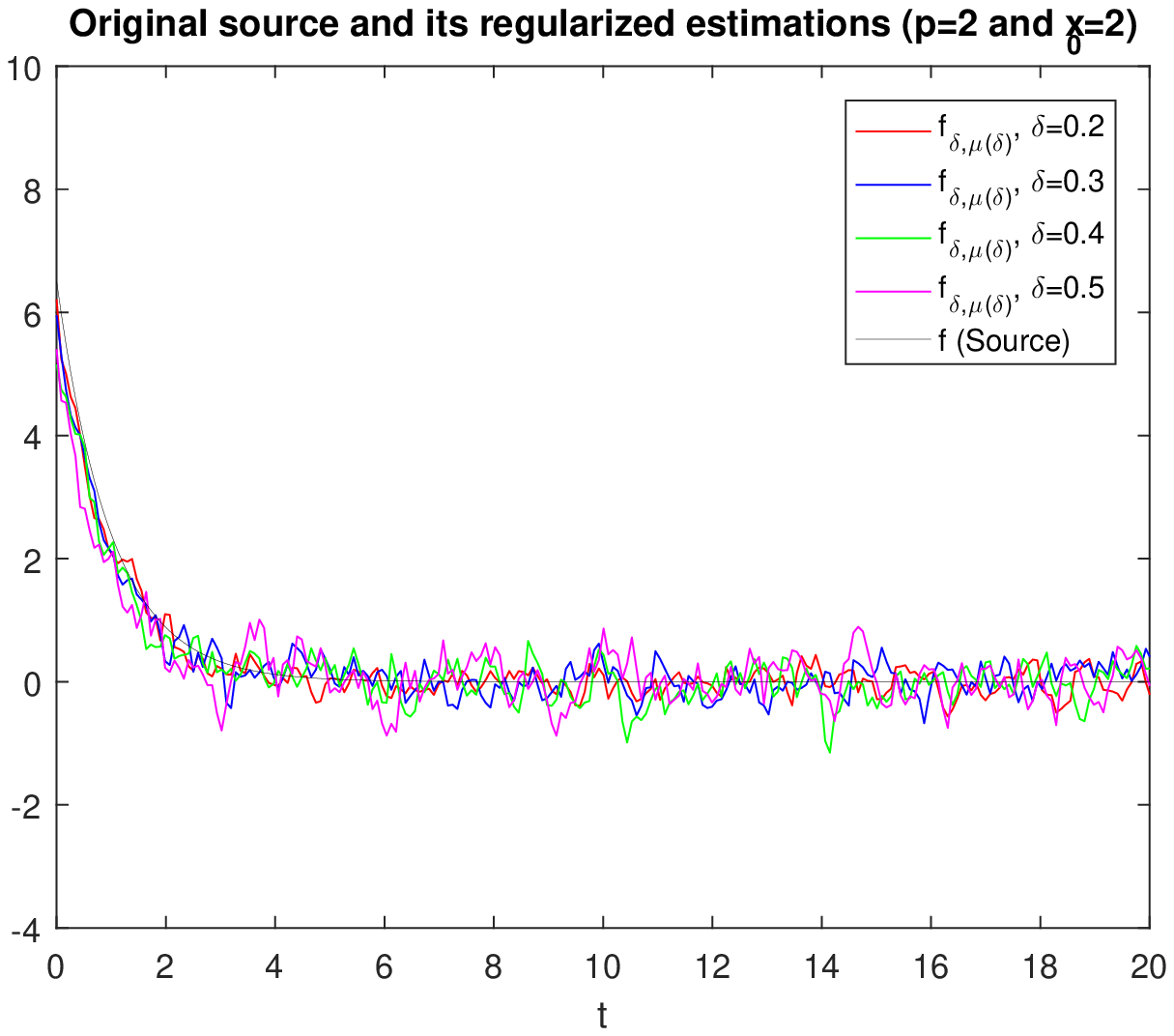}
\includegraphics[scale=0.55]{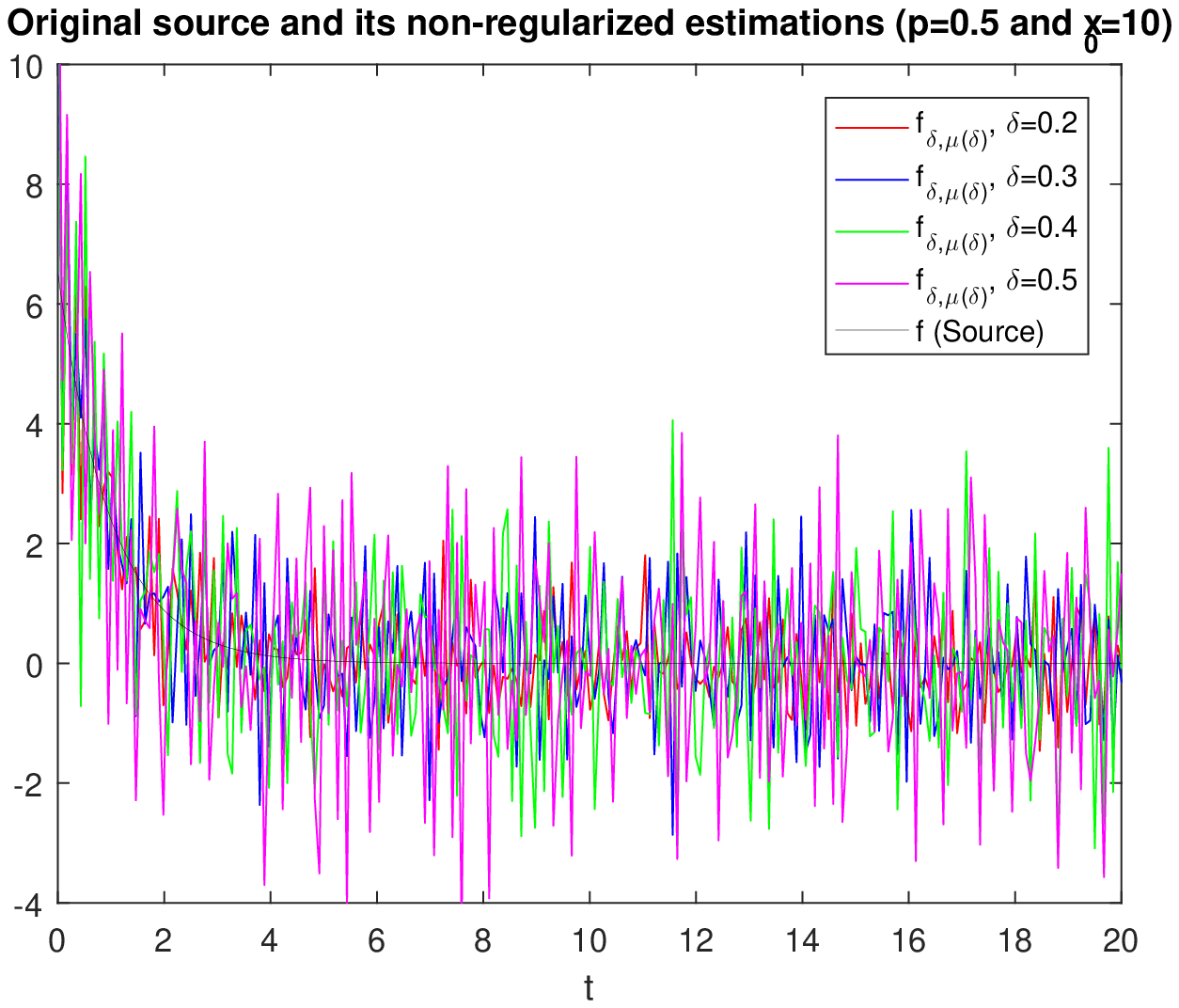}\hspace{-.5cm}
\includegraphics[scale=0.55]{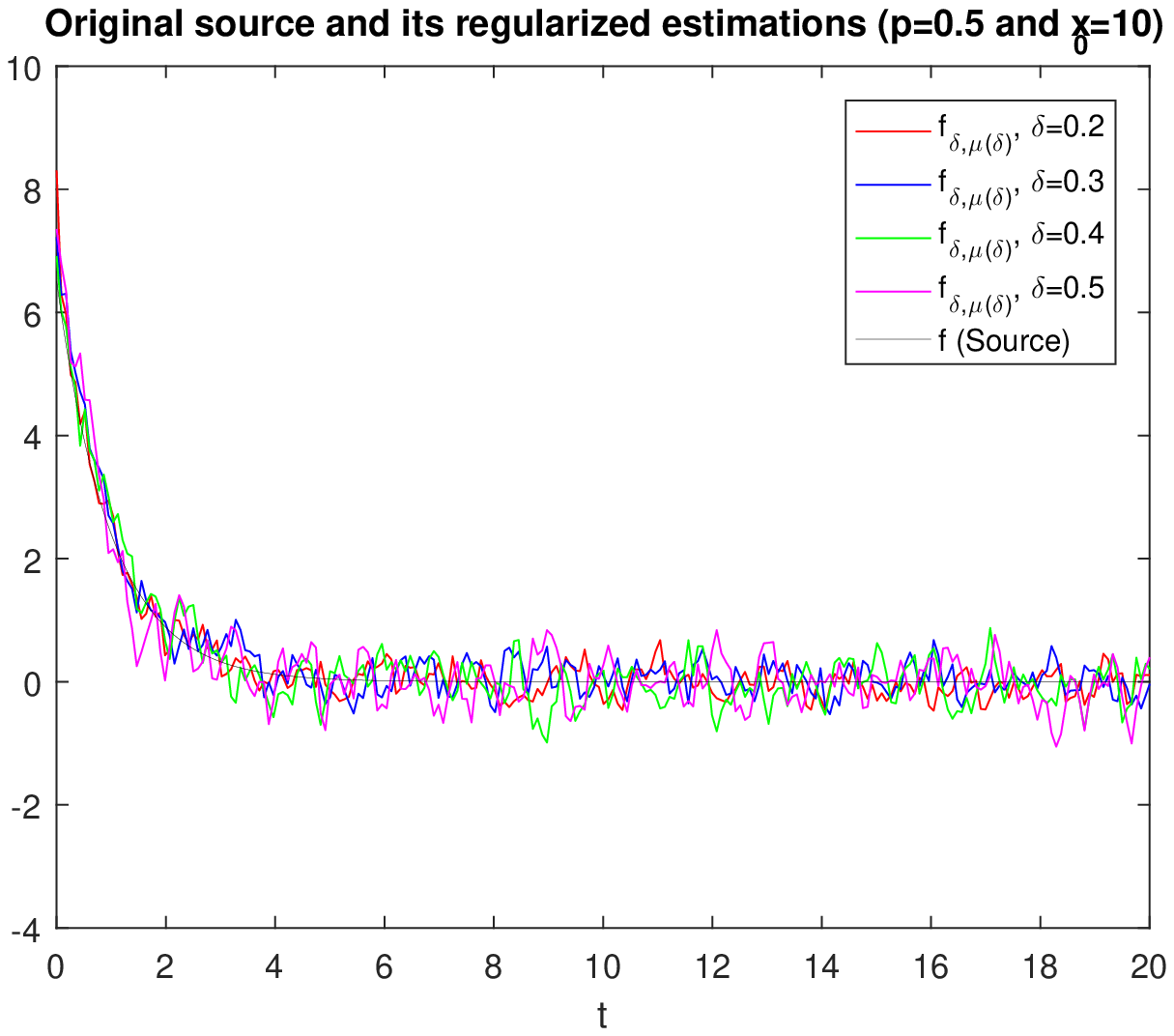}
\begin{center}
\caption{Example 1: Original Source and its estimations for the solution $f_\delta$ (left)  and the regularized one $f_{\delta,\mu}$ (right) for the same set of  noise levels.} 
\label{Exponential}
\end{center}
\end{figure}

\setcounter{theorem}{0}
We consider functions $u$ and $f$ that satisfy a given transport equations on  an interval $[0,T]$ and define a uniform partition ${\cal P}$ on that interval.  
We set the value of the regularization parameter $\mu =\delta^{\frac{2}{p+2}}$  where $\delta$ denotes the noise levels and calculate the approximated solution $f_{\delta ,\mu}(t)$  given by \eqref{ffreg} from simulated data $\{ u(x_0,t_j)=y(tj), t_j \in {\cal P} , \, j=0,...,n\}$. 
%
%
\begin{example}
\label{ej1}
Consider the inverse source problem defined in \eqref{transpeqn} with modeling parameter values  $\alpha^2=0.01, \beta=0.5, \nu= 1.51$. Hence, the problem is to determine $f$ in
$$\begin{cases} 
u_t(x,t)=0.01 u_{xx}(x,t)-0.5 u_x(x,t)-1.51 u(x,t)+f(t), \qquad & x \in \R, \, t>0 \\
u(0,t)=0, \qquad & t \in \R \\
u(x_0,t)=y(t), \qquad & t \in \R. \, \\ 
\end{cases}$$
from measured data $y_\delta(t)$ satisfying 
$\| y_\delta(t) -y(t)\|_{L^2(\R)} \leq \delta$ for given noise levels $\delta < 1$.
\end{example}
We consider the source $f$  given by 
\begin{equation}
\label{fexp}
f(t)=
\left\{ \begin{matrix}
6.51 e^{-t},\qquad & 20>t\geq 0, \\
0,\qquad &\, t\ge 20.
\end{matrix}
\right.
\end{equation}
\begin{example}
\label{ej2}
Consider the inverse source problem defined in \eqref{transpeqn} with modeling parameter values  $\alpha^2=0.1, \beta=0.9, \nu= 1$. Hence, the equation is given by
$$\begin{cases} 
u_t(x,t)=0.1 u_{xx}(x,t)-0.9 u_x(x,t)- u(x,t)+f(t), \qquad & x \in \R, \, t>0 \\
u(0,t)=0, \qquad & t \in \R \\
u(x_0,t)=y(t), \qquad & t \in \R.
\end{cases}$$
\end{example}
We numerically study the behavior of the regularization method for source $f$  given by 
\begin{equation}
\label{f2}
f(t)=
\begin{cases}
0, \qquad & 2>t\geq 0, \\
2, \qquad & 4>t\geq 2, \\
-2, \qquad & 6>t\geq 4, \\
1, \qquad & 8>t\geq 6, \\
-1, \qquad & 10>t\geq 8, \\
1/2, \qquad & 12>t\geq 10, \\
-1/2, \qquad & 14>t\geq 12, \\
1/4, \qquad & 16>t\geq 14, \\
-1/4, \qquad & 18>t\geq 16, \\
0, \qquad &  t\geq 18.
\end{cases}
\end{equation}
Figure \ref{Exponential} shows the original source $f$ along with the non-regularized solution $f_\delta = Ty_\delta$ (left side) and the regularized one $f_{\delta}$ given by \eqref{ffreg} (right side).
The plots on the top correspond to $p=2$ and the data measured at $x_0=2$ and the one below correspond to $p=0,5$ and $x_0=10$. Also, high level noise were considered, we take $\delta \in \{0.2, 0.3, 0.4, 0.5\}$. 
Table \ref{tabej1} contains the absolute errors for the non-regularized and the regularized source solutions when considering $p=2$ and measuring position $x_0=2$.\\

\begin{figure}[h!]
\begin{center}
\includegraphics[width=0.5\textwidth]{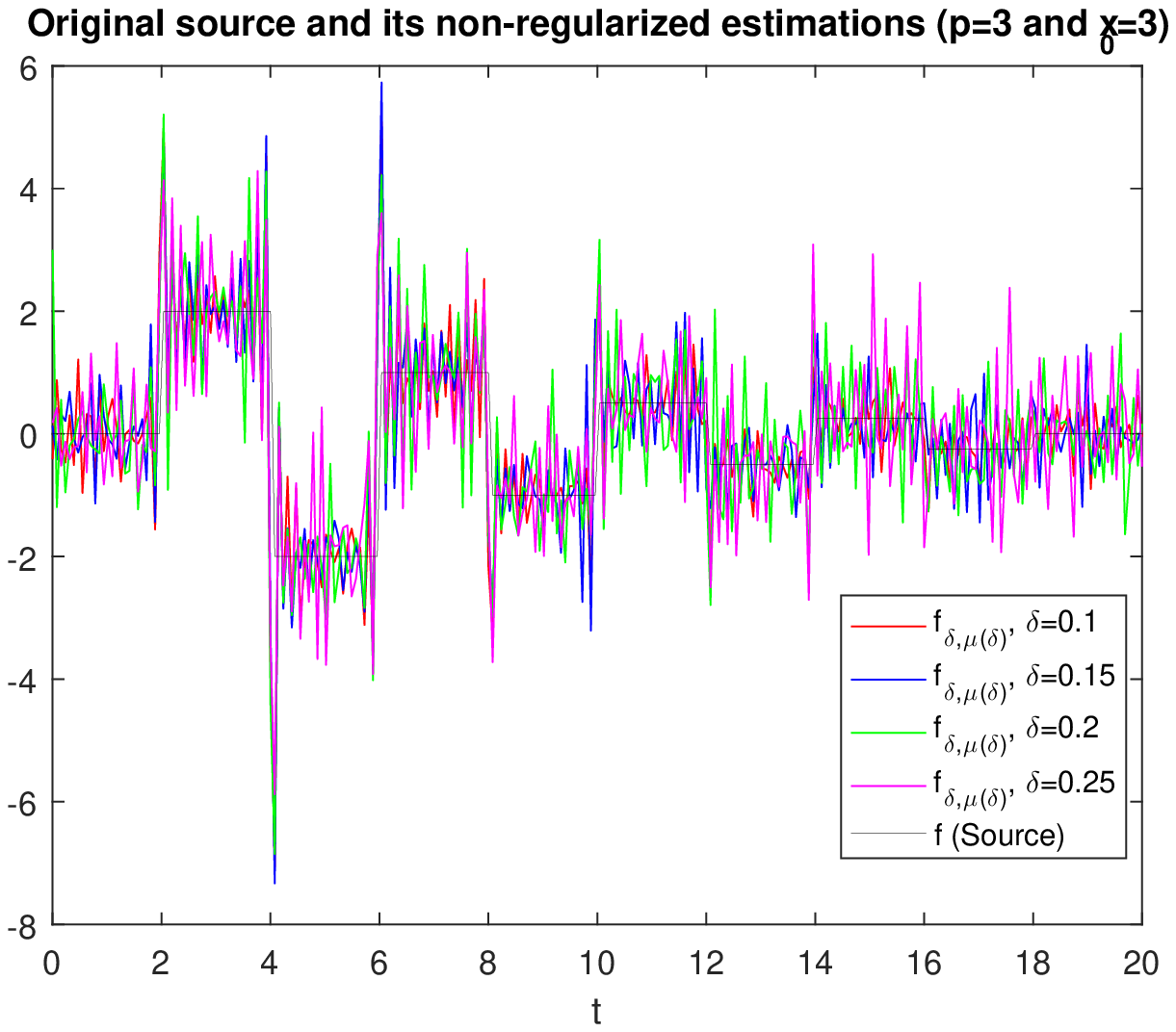}\hspace{-.25cm}
\includegraphics[width=0.5\textwidth]{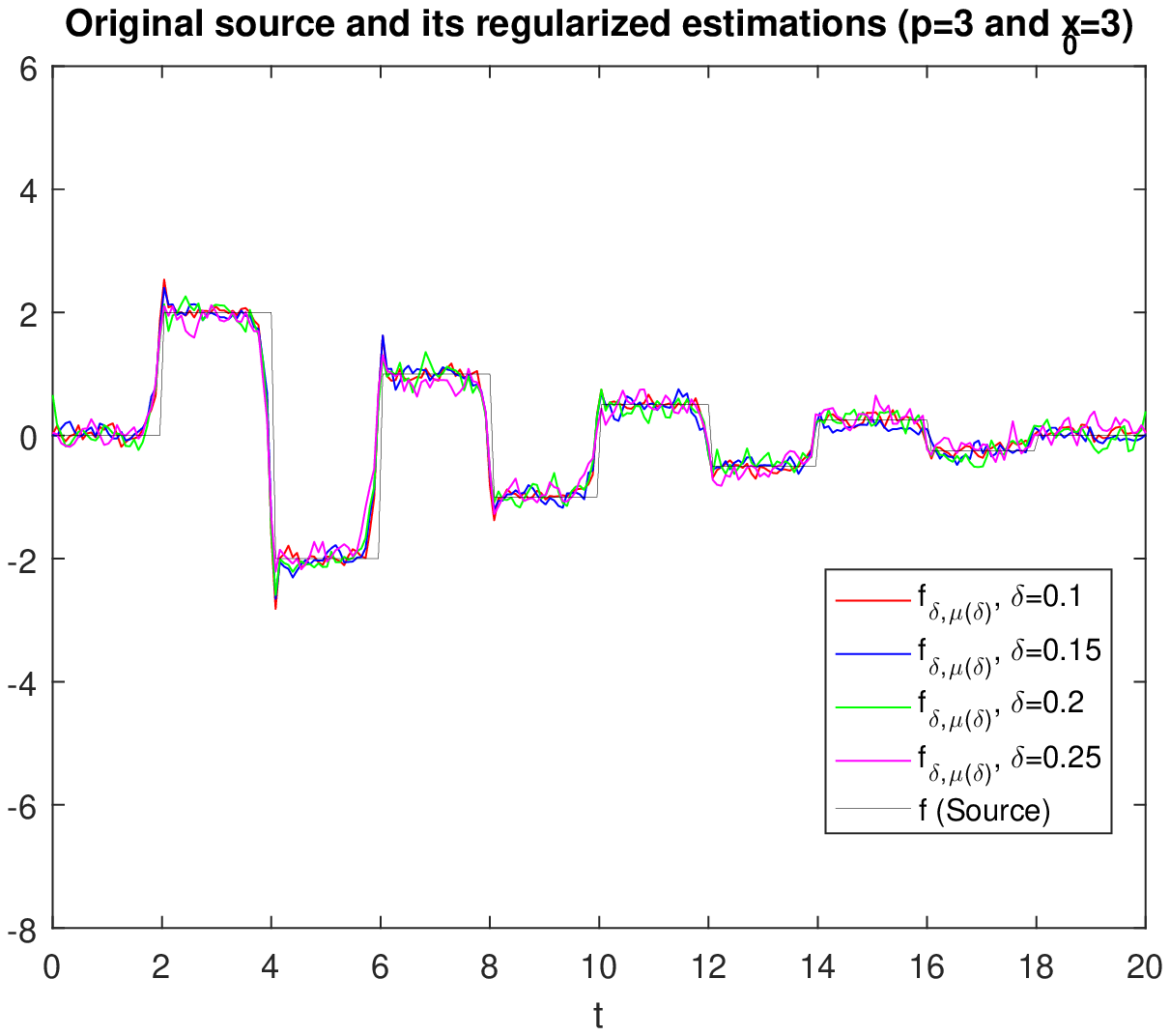}
\includegraphics[width=0.5\textwidth]{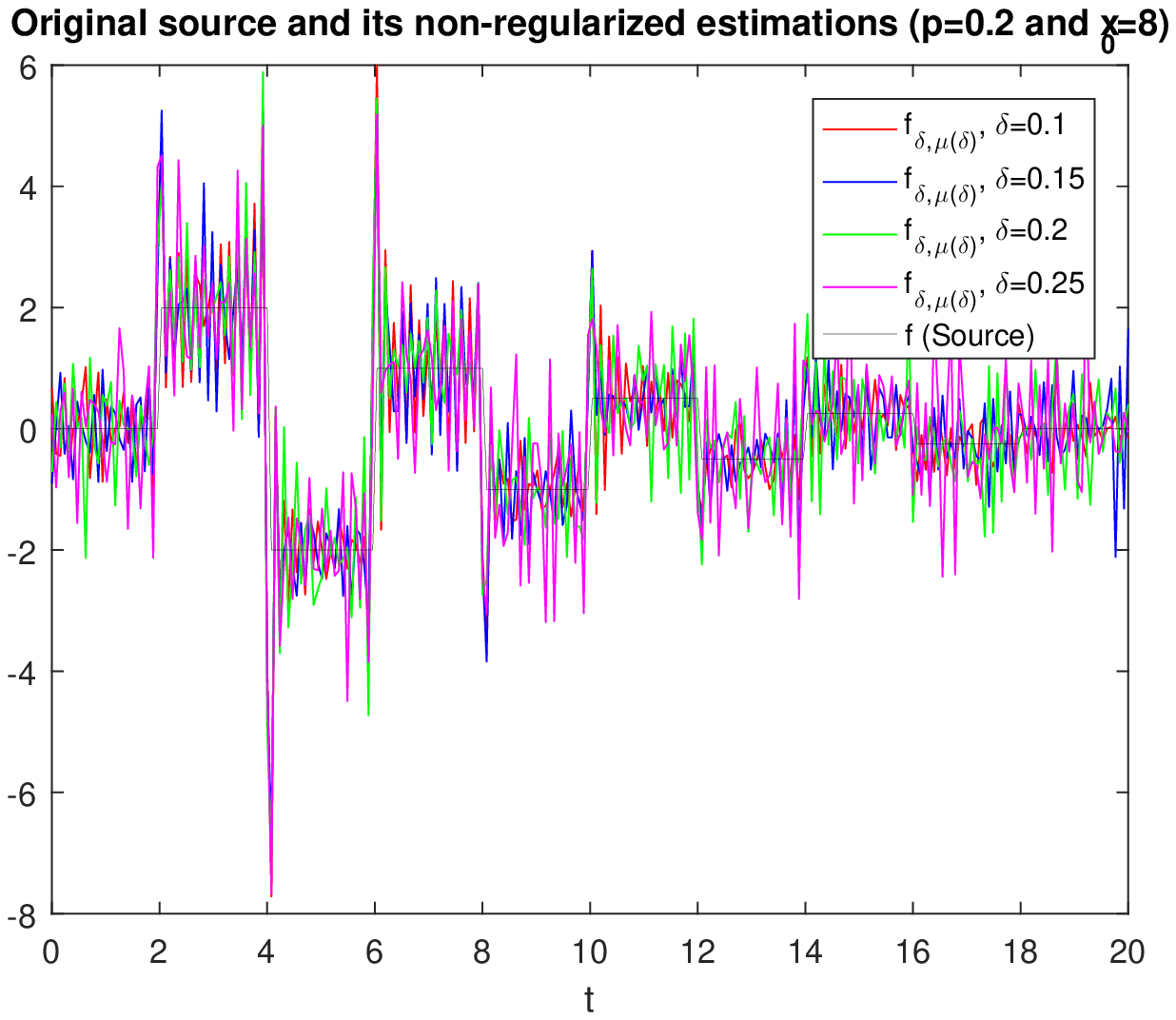}\hspace{-.5cm}
\includegraphics[width=0.5\textwidth]{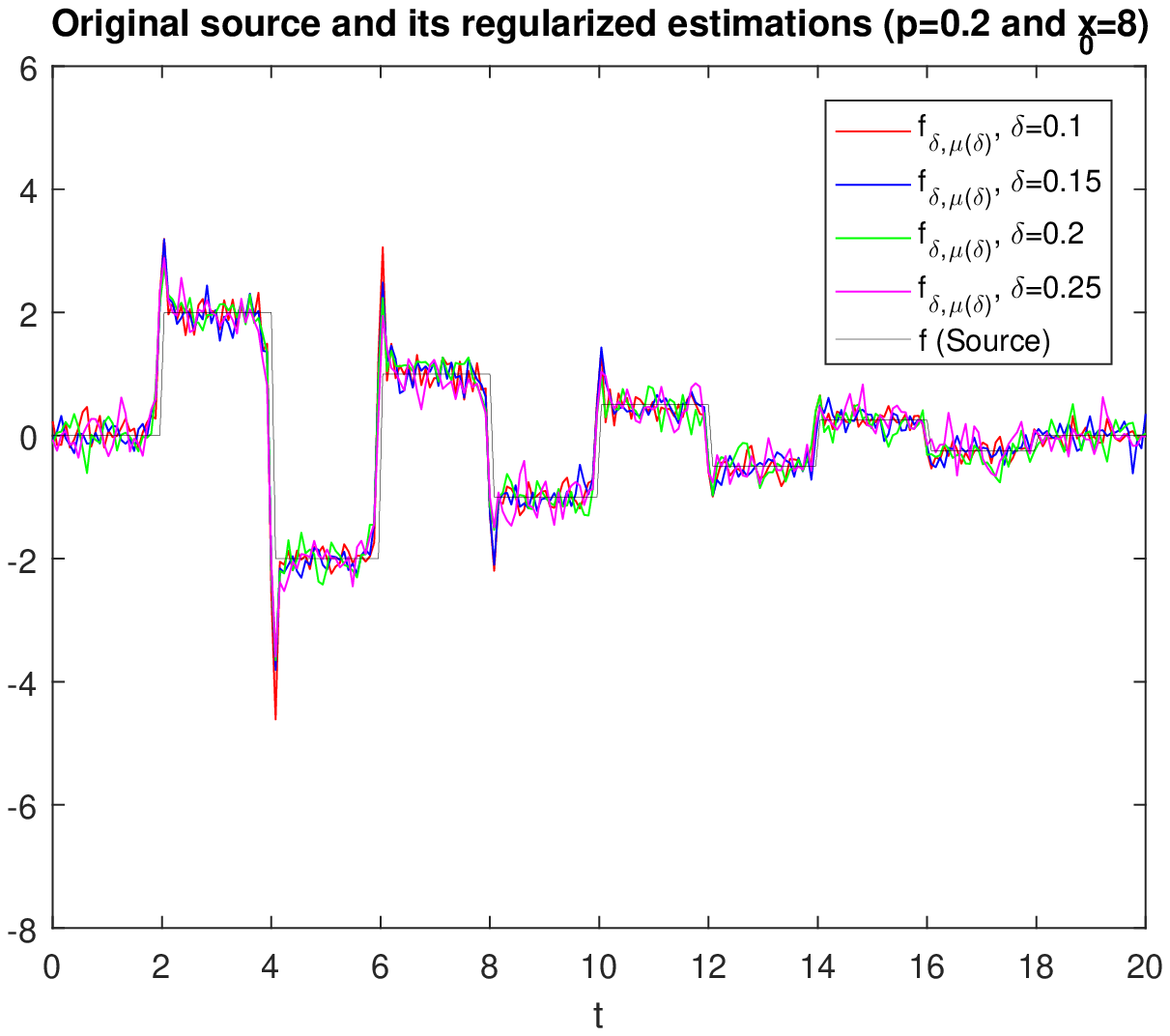}
\caption{Example 2: Original Source and its estimations $f_\delta$ (left) and the regularized solution $f_{\delta,\mu}$ (right) for the same set of noise levels.} 
\label{Ondacuadrada}
\end{center}
\end{figure}

Figura \ref{Ondacuadrada} shows the original source $f$ along with  $f_{\delta}$ (left side) and the regularized solution  $f_{\delta,\mu}$  (right side).
The plots on the top correspond to $p=3$ and the data measured at $x_0=3$ and the ones below correspond to $p=0.2$ and $x_0=8$. For this example, the level noises are  $\delta \in \{0.1, 0.15 0.2, 0.25\}$. 
Table \ref{tabej2} contains the absolute errors for the non-regularized and the regularized source solutions when considering $p=3$ and measuring position $x_0=3$.

\begin{table}[h!]
\begin{tabular}{|c|c|c|c|c|c|c|c|c|c|c|}
\hline $\delta$ & 0,01   & 0,02  & 0,03  & 0,04  & 0,05  & 0,06  & 0,07  & 0,08  & 0,09  & 0,1 \\
\hline
\hline
 $\|  f-f_{\delta} \|$ & 9,54 &	9,59	& 9,62	& 9,75	& 10,09	& 10,17	& 10,283	& 10,90	& 11,85	& 12,19 \\
\hline
 $\|  f-f_{\delta,\mu} \|$ & 5,282 &	5,305	& 5,307	& 5,309	& 5,310	& 5,315	& 5,318	& 5,320	& 5,326	& 5,369 \\
\hline
 {\tiny Theoretical Bound}  & 6,47	& 6,75 &	6,98 &	7,13 &	7,46	& 7,84 &	8,11 &	8,30 &	9,41 &	10,46\\
\hline
\end{tabular}
\caption{Example 1: Absolute non-regularized and regularized estimation errors  and theoretical bound for example \ref{ej1} assuming $p=2$ and data measured at $x_0=2$.}
\label{tabej1}
\end{table}  
\begin{table}[h!]
\centering
\begin{tabular}{|c|c|c|c|c|c|c|c|c|c|c|}
\hline $\delta$ & 0,01   & 0,02  & 0,03  & 0,04  & 0,05  & 0,06  & 0,07  & 0,08  & 0,09  & 0,1 \\
\hline
\hline
 $\|  f-f_{\delta} \|$ & 13,64	& 13,72	& 13,82	& 13,88	& 13,97	& 14,01	& 14,23	& 14,61	& 14,72	& 15,43  \\
\hline
 $\|  f-f_{\delta,\mu} \|$ & 2,21	& 2,85	& 3,23	& 3,97	& 4,31	& 4,78	& 4,92	& 5,03	& 5,49	& 5,88  \\
\hline
{\tiny Theoretical Bound} & 2,37	& 3,13	& 3,68	& 4,13	& 4,52	& 4,86	& 5,17	& 5,46	& 5,72	& 5,97	\\
\hline
\end{tabular}
\caption{Example 2: Absolute non-regularized and regularized estimation errors and theoretical bound for example \ref{ej2} assuming $p=3$ and data measured at $x_0=3$.}
\label{tabej2}
\end{table}  

\section{Conclusions}
\label{conclu}

 
We consider the inverse source problem for a 1D transport equation.
We define a regularization family of operators  to deal with the ill-posedness of the problem  by compensation the instability factor in the inverse operator. We proposed a regularization parameter choice rule based on assumption of the noise level in data and the smoothness of the source to be identify. We prove that for the parameter choice rule proposed here, the method is stable and a H\"older type  bound for the regularization error is obtained.
The numerical examples show an improvement in the regularized solution with the respect to the one obtained when no regularization is applied.
Numerical examples  show good estimates for the source at different noise levels, in this work  we have included  few cases where the sources belong to different Hilbert spaces to illustrate the performance of this method.

\end{document}